\documentclass{article}

 \usepackage{tikz}
\usepackage{amsthm,amssymb,amsmath}
\usepackage{epsfig}
\usepackage{color}
\usepackage{algorithm2e}
\usepackage{algpseudocode}
\usepackage{subcaption}
\usepackage{ mathrsfs }
\newtheorem{axiom}{Axiom}

 \newtheorem{thm}{Theorem}

\newtheorem{example}{Example}
\theoremstyle{definition}
\newtheorem{defn}[thm]{Definition}
\theoremstyle{remark}
\newtheorem{rem}{Remark}

\title{Optimal Reinsurance: A Ruin-Related Uncertain Programming Approach}
\author {Wrya Vakili  and Alireza Ghaffari-Hadigheh\\
Dept. of Applied Math. \\
 Azarbaijan Shahid Madani University, Tabriz, Iran\\
{\tt vrichvri@gmail.com}\\
{\tt hadigheha@azaruniv.ac.ir}
}
\date{}
\begin{document}
\maketitle

\begin{abstract}
We investigate the role of reinsurance in maximizing the wealth of an insurance company. We use Liu's uncertainty theory (B. Liu, 2007) for the problem modeling and follow-up computations. The uncertainty measure of ruin for the insurance company is considered as the optimization criterion.  Since calculating the ruin index is very difficult, we introduce a simple computational method to identify the uncertain measure of ruin for an insurance company. Finally, a generalized model is presented, granting the model be more practical.
\end{abstract}

\section{Introduction}
Risk is an inevitable feature of a business. Analyzing the behavior of involved risk and it's consequences have always been subjected to different theories that deal with unpredictable environments. Risk theory investigates the effect of potential deviation of the outcome from their expected value,   prevent it from occurring, or being prepared to take action when it is inescapable. It mostly practices the modeling of cash flow as a surplus process, e.g., the business's wealth during the activity period \cite{bowers1987actuarial}.

Using insurance, people aim to hedge their property against diverse sorts of risks they might encounter. Insurance companies take their customers' risks in exchange for a specified amount of money referred to as the premium. The collective risk model, also called the Cramer-Lundberg model, was initiated by Lundberg and latter developed by Cramer \cite{asmussen2010ruin}. It is usually counted as a standard configuration of an insurance company's risk process concerning the risks it takes by selling insurance underwriting.

 In the non-life insurance industry, consider an insurance company that sells many finite fire insurance contracts. To compete against other companies, it may offer a lower premium, leading to enduring a high risk in each contract. Notice that in such contracts, the chance of damage of property by the fire is significantly low. However, the policyholder's loss in an incident is extremely higher than the premium amount. Exercising a mass of such a risk, the company would prefer to hedge itself by purchasing reinsurance to transform some of the risks to another company. Here, we model such a situation as an optimization problem at which the optimal decision is the amount of risk it should share.

Consider  this classical model as appeared in most of the standard actuarial literature \cite{asmussen2010ruin}
\begin{equation}
U_t=u+ct -R_t,
\end{equation}
where  $u$ is the initial capital of the company, $c > 0$ is a constant that represents insurance premium indicating that company receives $c$ units of money per unit time. The process $R_t =\sum_{i=1}^{N_t}\eta_i$ is so-called compound process, where $N_t$ is the process of counting the number of claims that the company might have to pay during $(0,t]$. Further,  $\eta_i$'s are random variables denoting the severity of claims. Finally, $U_t$ is the wealth of the company regarding its contracts over $(0,t]$.

Observe that the interval $(0,t]$ can be partitioned into separated points  $s_1,s_2,\ldots $ in which losses take place.  We  refer to  $\xi_i=s_{i+1}-s_i$  as interarrival times. In the  Lundberg-Cramer model, these random variables are iid with exponential distribution. It is also assumed that $N_t$ and $\eta_i$'s are independent, and $\eta_i$'s are iid random variables.
Anderson et al.\cite{andersen1957collective} generalized the classic model by letting $N_t$ to be any renewal process, not just the Poisson process. In these models, the only non-deterministic component is $R_t$, while calculating its distribution is a very costly task, and most of the time, it is almost impossible to do. However,  an acceptable distribution can be approximated by discretizing the variables to reduce them to be more tractable.

Ruin is an essential and primary concept of mathematical approaches in the risk management of insurance activities. It is defined as the company's surplus being negative in some future, concerning active contracts of some policies \cite{asmussen2010ruin}. Observe that ruin happens when the collected premiums don't satisfy the losses claimed by policyholders.  Just by restricting the model to some simplified considerations, it is possible to obtain some results.  For example, Panjer \cite{panjer1981recursive} proved that when severity only takes positive integers, and $N_t$ is of some special cases,  calculation of the distribution of $R_t$ and consequently the ruin, is possible through some recursive algorithm.
In another study,  diffusion approximation was proposed  \cite{iglehart1969diffusion} and followed by Grandell \cite{grandell1977class}, as another approach to make the computation of ruin possible. Here, the main idea is using the central limit theorem in favor of approximating the risk process by a Wiener process. Another approximation was proposed by De Vydler et al. \cite{de1}, where the first three moments of the $R_t$ were used to mimic its behavior and find a recursive algorithm for calculating the probability distribution of the ruin.
Thus, direct calculating of the probability of ruin is far from being a practical approach.

  Some generalized versions of classical models were developed, where more details of the real-world problem are taken into account. Dufresne et al. \cite{dufresne1991risk} used a diffusion process to perturb the risk process, allowing more details to the model, such as the effect of interest rate or dependency between the number and severity of the claims. More assumptions, such as taxes and dividends, have been added to models enriching the literature and making the models more practical. We refer to  \cite{albrecher2007lundberg,yang2005optimal}, as well as to the masterpiece textbooks in the field \cite{asmussen2010ruin,embrechts2013modelling} for more details.

Reinsurance is regarded as a framework of protecting the insurance companies against potential catastrophic high losses \cite{albrecher2017reinsurance}. As an insurance company charges its policyholders with premiums, the reinsurer also defines a premium adapted to the amount of preassumed uncertain compensation cover that offers the cedent party. Reinsurance's role is to limit the primary insurance liability with a predefined risk and increase its risk-taking capacity when the likely loss is predicted to be overwhelming.

There are two basic forms of reinsurance contracts; proportional and non-proportional. In the former, both parties admit a predetermined share of all premiums and their potential losses.  In the latter,  reinsurance's liability is not determined already but depends on the number or severity of the claims incurred before the reinsurance contract's expiration~\cite{albrecher2017reinsurance}.

One of the most important concerns about the reinsurance contracts is determining optimal retention, i.e., the optimum amount of the risk that primary insurance holds. Since the pioneering works by Borch \cite{borch1960safety}, Kahn \cite{kahn1961some}, and Arrow \cite{arrow1978uncertainty}, there have been many contributions to the subject. Followed by these fundamental works, other researchers added a great deal of knowledge to the literature, e.g.,   \cite{bowers1987actuarial,buhlmann2007mathematical,gerber1979introduction,van1989optimal,waters1979excess,waters1983some}.

The special case of proportional reinsurance was the subject of some other researches \cite{cai2007optimal,cai2008optimal,centeno1985combining,kaluszka2001optimal}. Common objectives  include minimizing the probability of ruin \cite{bai2013optimal,chen2010optimal,promislow2005unifying}, maximizing the expected utility of the terminal wealth  \cite{bai2010optimal,gerber2019constraint,liu2009optimal,liang2014optimal}, and mean variance portfolio optimization \cite{bauerle2005benchmark,kaluszka2001optimal}. For a comprehensive review on optimal reinsurance approaches see \cite{de2009optimal}.

Probability theory based models are assumed to be practical and useful, not to be considered the perfect and the only way to model certain problems. As noted by some people \cite{liu2012there,zadeh1995discussion,zadeh2008there}, probability theory has some shortcomings in dealing with those problems with no sufficient data, or basic information is based on some experts' opinions.

Liu \cite{liu2007uncertainty} introduced the uncertainty theory to deal with problems related to belief degree in expert-based systems. Analogous to the stochastic point of view, Liu \cite{liu2008fuzzy} introduced the concept of renewal and renewal reward process for initiated systems based on sudden arrival events. To deal with the ruin theory problems, the concept of ruin index has been proposed in the uncertainty theory framework as its substitution in probability theory \cite{liu2013extreme}.  Li et al. \cite{li2013uncertain} proposed a  premium calculation principle in an uncertain environment. Later, Yao et al. \cite{yao2015modified,yao2016ruin} modified the results on ruin index in some extensions and generalized over different assumptions.\medskip

Uncertainty theory has two crucial features that make it a reliable tool dealing with expert-based models. Firstly, it enjoys a well-defined axiomatic structure similar to probability theory. Secondly,  this theory is designed in a way to be tractable for subjective modeling in the fields like actuarial science that distributions are highly based on belief and not data. Here we investigate the impact of proportional reinsurance on an insurance company's surplus process and its survival that is defined as leaving the uncertain measure of ruin (UMR) under a fair bound. We assume a quota share reinsurance contract, through which the primary insurance cedes premiums and losses with an unknown cession rate to the reinsurer. We need to consider a safety loading for both reinsurance and insurance premiums so that the UMR does not occur in an almost certain way in the long run.  We study the optimal retention of proportional reinsurance under ruin-related optimization restrictions. The model determines parameters that identify the best possible scenario to get involved in a reinsurance contract. Based on a proportional type of reinsurance contract, we only study the problem from an insurance point of view, assuming that reinsurers are willing and capable of taking the whole committed risk if required.

The rest of the paper is organized as follows. Since uncertainty theory is relatively new, in the next section we present some basics on the theory which will be used in the main parts of the paper. Section 3 includes three subsections. First, we introduce a new way of calculating the uncertain measure of ruin followed by a detailed example. We then discuss optimal reinsurance as an uncertain optimization problem and describe the results through an example in non-life insurance. The section will be ended with a generalization of the proposed model. The conclusion suggests some possible extensions that might make the model more reliable.

\section{Preliminaries } \label{prelim}
 Let $\Gamma $ be a nonempty set and  $\mathscr{L}$ be a $\sigma$- algebra over $\Gamma $. Each element $\Lambda \in \mathscr{L}$ is called an event. A set function $\mathscr{M}$ from $\mathscr{L}$ to $[0,1]$ that satisfies the following axioms is called an uncertain measure \cite{liu2007uncertainty}:

\begin{axiom} (Normality Axiom) $\mathscr{M}\lbrace \Gamma \rbrace  =1$.
\end{axiom}
\begin{axiom} (Duality Axiom)  $\mathscr{M}\lbrace \Lambda \rbrace+\mathscr{M}\lbrace \Lambda^{c} \rbrace =1$ for any event $\Lambda$.
\end{axiom}
\begin{axiom} (Subadditivity Axiom) For every countable sequence of events $\Lambda_1, \Lambda_2, \ldots$,
\begin{equation}
\mathscr{M}\lbrace \bigcup_{i=1}^{\infty} \Lambda_{i}\rbrace  \leq \sum_{i=1}^{\infty}\mathscr{M}\lbrace \Lambda_{i}\rbrace .
\end{equation}
\end{axiom}
The triplet $(\Gamma,\mathscr{L}, \mathscr{M})$ is called an uncertainty space.\\
The Product axiom, distinguishing the probability theory from the uncertainty theory, is defined as follows~\cite{Liusome}.

\begin{axiom} (Product Axiom) Let $(\Gamma_k,\mathscr{L}_k,\mathscr{M}_k)$ be uncertainty spaces for $k=1,2,\ldots$. The product uncertain measure $\mathscr{M}$ is a measure on product $\sigma$-algebra $\mathscr{L}_1\times \mathscr{L}_2\times\cdots \times \mathscr{L}_n $ satisfying
\begin{equation}
\mathscr{M}\lbrace \prod_{k=1}^{\infty}\Lambda_k \rbrace=\bigwedge_{k=1}^{\infty}\mathscr{M}_{k}\lbrace\Lambda_k \rbrace.
\end{equation}
\end{axiom}
 Uncertain variable has been  defined  to ease the quantitative modeling of phenomena  in uncertainty theory~\cite{liu2007uncertainty}. It is a measurable function from an uncertain space $(\Gamma,\mathscr{L},\mathscr{M})$  to the set of real numbers, in which for any Borel set $B$, the set $\lbrace \xi \in B\rbrace=\lbrace \gamma\in\Gamma \vert \xi(\gamma)\in B\rbrace$ is an event.
  Uncertainty distribution of an uncertain variable $\xi$  is defined as \cite{liu2007uncertainty}
 \begin{equation}
\Phi(x)=\mathscr{M}\lbrace \xi\leq x \rbrace, \qquad \forall x\in \Re.
\end{equation}
There are several  uncertain variables in the literature of the uncertainty theory. The simplest one is the linear uncertain variable with the following distribution,

\begin{align}
 \Phi(x)= \begin{cases}
      0 &~{\rm if} ~ x\leq a \\
      \displaystyle\frac{x-a}{b-a} &~{\rm if} ~ a\leq x\leq b \\
      0 & ~{\rm if} ~ x\geq b,
   \end{cases}
\end{align}

denoted by $\mathscr{L}(a,b)$ where $a$ and $b$ are real numbers with $a < b$.

An uncertain variable $\xi$ is called normal and denoted  by $\mathscr{N}(e,\sigma)$ \cite{liu2007uncertainty}, if it has the normal uncertain distribution
\begin{equation}
\Phi(x)=\left(1+exp\left(\frac{\pi(e-x)}{\sqrt{3}\sigma}\right)\right)^{-1},x\in \Re,
\end{equation}
 where $e$ and $\sigma$ are real numbers with $\sigma > 0$. Moreover, an uncertain variable $\xi$ is called lognormal if $\ln \xi$ is a normal uncertain variable. Its distribution reads as
 \begin{equation}
 \Phi(x)=\left(1+exp\left(\frac{\pi(e-\ln x)}{\sqrt{3}\sigma}\right)\right)^{-1}, x\geq 0.
 \end{equation}
 Uncertain variables $\xi_1,\xi_2,\dots, \xi_n$ are  independent if \cite{Liusome}
\begin{equation}
\mathscr{M}\left\lbrace \bigcap_{i=1}^{n}( \xi_i \in B_i ) \right\rbrace = \bigwedge_{i=1}^{n}\mathscr{M}\left\lbrace \xi_i \in B_i \right\rbrace,
\end{equation}
for arbitrary Borel sets $B_1,B_2,\dots, B_n$.
An uncertain distribution $\Phi(x)$ is called regular if it is a continuous and strictly increasing function w.r.t.  $x$ for all $0<\Phi(x)<1 $, and
\begin{equation}
\lim_{x\rightarrow -\infty}\Phi(x)=0 ,\qquad  \lim_{x\rightarrow +\infty}\Phi(x)=1.
\end{equation}
Let $\xi_1, \xi_2,\ldots,\xi_n$ be independent uncertain variables with regular uncertainty distributions $\Phi_1,\Phi_2 ,\ldots,\Phi_n$, respectively. If $f(\xi_1, \xi_2,\ldots,\xi_n)$ is strictly increasing w.r.t. $\xi_1, \xi_2,\ldots,\xi_m$ and strictly decreasing w.r.t.  $\xi_{m+1}, \xi_{m+2},\ldots,\xi_n$, then $f$ is an uncertain variable with the uncertainty distribution \cite[Theorem 1.26]{Version10}
\begin{equation}
\Psi(x)=\sup_{f(x_1,x_2,\ldots,x_n)=x}\left( \min_{1\leq i\leq m }\Phi_i(x_i)\wedge \min_{m+1\leq i\leq n}(1-\Phi(x_i)) \right).
\end{equation}
\begin{defn}
Let $\xi$ be an uncertain variable with regular uncertainty distribution $\Phi(x)$. The inverse function $\Phi^{-1}(\alpha)$ is called the inverse uncertainty distribution of $\xi$.
\end{defn}
\begin{thm}\label{Th1}
 If $f(\xi_1, \xi_2,\ldots,\xi_n)$ is strictly increasing w.r.t. $\xi_1, \xi_2,\ldots,\xi_m$ and strictly decreasing w.r.t.  $\xi_{m+1}, \xi_{m+2},\ldots,\xi_n$, then $f$ is an uncertain variable with the inverse uncertainty distribution
\begin{equation}\label{Eq1}
\Psi^{-1}(\alpha)=f(\Phi^{-1}(\alpha),\ldots, \Phi^{-1}(\alpha), \Phi_{m+1}^{-1}(1-\alpha),\ldots, \Phi_n^{-1}(1-\alpha)).
\end{equation}
\end{thm}
\begin{thm}\label{Theorem11} \cite{liu2010uncertain}
Let $\xi_1,\xi_2,\ldots ,\xi_n$ be independent uncertain variables with regular uncertainty distributions $\Phi_1,\ldots ,\Phi_n$, respectively. If $f(\xi_1,\ldots ,\xi_n)$ is strictly increasing with respect to $\xi_1,\xi_2,\ldots ,\xi_m$ and strictly decreasing with respect to $\xi_{m+1},\xi_{m+2},\ldots ,\xi_n$, then
\begin{equation}
\mathscr{M}\left\lbrace f(\xi_1,\xi_2,\ldots ,\xi_n)\leq 0  \right\rbrace ,
\end{equation}
is the root $\alpha$ of the equation
\begin{equation}
f(\Phi_1^{-1}(\alpha),\ldots ,\Phi_m^{-1}(\alpha),\Phi_{m+1}^{-1}(1-\alpha),\ldots ,\Phi_{n}^{-1}(1-\alpha))=0.
\end{equation}
\end{thm}
\noindent
If $f(\Phi_1^{-1}(\alpha),\ldots ,\Phi_m^{-1}(\alpha),\Phi_{m+1}^{-1}(1-\alpha),\ldots ,\Phi_{n}^{-1}(1-\alpha))< 0,$
for all $\alpha$, then we set the root $\alpha =1$; and if
$f(\Phi_1^{-1}(\alpha),\ldots ,\Phi_m^{-1}(\alpha),\Phi_{m+1}^{-1}(1-\alpha),\ldots ,\Phi_{n}^{-1}(1-\alpha))>0$
for all $\alpha$, then we set the root $\alpha =0$.

The expected value of an uncertain variable $\xi$ is defined as \cite{liu2007uncertainty}
\begin{equation}
E[\xi] = \int_0^{+\infty} \mathscr{M}\{\xi> r\} dr - \int_{-\infty}^0  \mathscr{M}\{\xi\leq r\} dr,
\end{equation}

provided that at least one of the two integrals is finite. Having an uncertain distribution $\Phi(x)$ of an uncertain variable $\xi$, one can  calculate its expected value by \cite{liu2007uncertainty}
\begin{equation}
E[\xi]=\int_{0}^{\infty}(1-\Phi(x))dx-\int_{-\infty}^{0}\Phi(x)dx,
\end{equation}
provided that at least one of these integrals is finite. It is proved that for an uncertain variable $\xi$ with regular uncertain  distribution $\Phi(x)$, it holds
\begin{equation}
E[\xi]=\int_{0}^{1}\Phi^{-1}(\alpha)d\alpha.
\end{equation}
 Variance of an uncertain variable $\xi$ with finite expected value $E[\xi]$ is also defined as \cite{liu2007uncertainty}
\begin{equation}
V[\xi]=E\left[ (\xi -E[\xi])^2 \right].
\end{equation}
Let $\xi$ and $ \eta$ be independent uncertain variables with finite expected values. Then, for any real numbers $a$ and $b$ \cite{Version10}
\begin{equation}\label{Expected}
E\left[ a\xi + b\eta \right] =a E\left[ \xi \right] +b E\left[ \eta \right].
\end{equation}
\subsection*{Uncertain process}
Let $(\Gamma_k,\mathscr{L}_k,\mathscr{M}_k)$ be  uncertainty spaces for $k=1,2, \ldots$,  and  $T$ be a totally ordered  set (e.g. time). For the sake of simplicity, we use the term ``time'' for each member of this  set. An uncertain process is a function $X_t(\gamma)$ from $T\times (\Gamma_k,\mathscr{L}_k,\mathscr{M}_k) $ to the set of real numbers such that $\lbrace X_t \in B \rbrace$ is an event for any Borel set $B$ of real numbers at each time $t$. An uncertain process $X_t$ has independent increments if
 \begin{equation}
 X_{t_0}, X_{t_1}- X_{t_0},X_{t_2}- X_{t_1},\ldots,X_{t_k}- X_{t_{k-1}},
 \end{equation}
 are independent uncertain variables where $t_0$ is the initial time and $t_1 , t_2,\ldots,t_k$ are any times with $t_0<t_1< t_2< \cdots < t_k$.

 An uncertain process has stationary increments if its increments are independent identically distributed  (iid) uncertain variables whenever the time intervals have equal lengths. For uncertain process $X_t $, fixing the event $\omega$, the resulted function  $X_t(\omega)$ is called a sample path of $X_t$.
\begin{defn}\cite{Liusome}\label{Liuprocess}
An uncertain process $C(t)$ is called Liu process if\\
(a) $C(0)=0$  and almost all sample paths are Lipschitz continuous,\\
(b) $C(t)$ has stationary and independent increments,\\
(c) Every increment $C(t+s)-C(s)$ is a normal uncertain variable with expected value zero and variance $t^2$.
\end{defn}
\begin{defn} \cite{liu2008fuzzy}
Let $\xi_1,\xi_2,\ldots  $ be iid uncertain interarrival times. Define $S_0=0$ and $S_n = \xi_1+\xi_2+\cdots+\xi_n $ for $n \geq 1$. The uncertain process
\begin{equation}
N_t=\max_{n\geq 0}\lbrace n \vert S_n\leq t \rbrace
\end{equation}
is called an uncertain renewal process.
\end{defn}
\begin{thm} \cite{Version10}
Let $N_t$ be a renewal process with iid uncertain interarrival times  $\xi_1,\xi_2,\ldots  $. Further, let $\Phi$ denote their common uncertainty distribution. Then, $N_t$ has an uncertainty distribution
\begin{equation}
\Upsilon_t(x)=1-\Phi\left( \frac{t}{\lfloor x\rfloor+1}\right), \qquad \forall x\geq 0,
\end{equation}
where $\lfloor x\rfloor$ represents the floor function.
\end{thm}
\begin{thm} \cite{Version10}
Let $N_t$ be a renewal process with iid uncertain interarrival times $\xi_1,\xi_2,\ldots  $ Then the average renewal number
\begin{equation}
\frac{N_t}{t}\rightarrow \frac{1}{\xi_1}
\end{equation}
in the sense of convergence in distribution as $t\rightarrow \infty$.
\end{thm}
\begin{thm} \cite{Version10}
Let $N_t$ be a renewal process with iid uncertain interarrival times $\xi_1,\xi_2,\ldots  $ Then
\begin{equation}
\lim_{t\rightarrow \infty}\frac{E[N_t]}{t}=E\left[  \frac{1}{\xi_1}\right] .
\end{equation}
If $\Phi$ is regular, then
\begin{equation}
\lim_{t\rightarrow \infty}\frac{E[N_t]}{t}=\int_{0}^1 \frac{1}{\Phi^{-1}(\alpha)}d\alpha.
\end{equation}
\end{thm}
\begin{defn}\cite{Version10}
Let $\xi_1,\xi_2,\ldots  $  be iid uncertain interarrival times. Further, let  $\eta_1,\eta_2,\ldots $ be iid uncertain  rewards or costs ( losses in insurance case) associated with $i$-th interarrival times $\xi_i$ for $i=1,2,\cdots$Then
\begin{equation}
R_t=\sum_{i=1}^{N_t}\eta_i,
\end{equation}
is called a renewal reward process, where $N_t$ is the renewal process with uncertain interarrival times $\xi_1,\xi_2,\ldots  $
\end{defn}
\begin{thm}\cite{Version10}
Let $R_t$ be a renewal reward process with iid uncertain interarrival times $\xi_1,\xi_2,\ldots  $ and iid uncertain rewards  $\eta_1,\eta_2,\ldots $ Assume $(\xi_1,\xi_2,\ldots )$ and $(\eta_1,\eta_2,\ldots )$ are independent uncertain vectors, and those interarrival times and rewards have uncertainty distributions $\Phi$ and $\Psi$, respectively. Then, $R_t$ has  uncertainty distribution
\begin{equation}
\Upsilon_t(x)=\max_{k\geq 0}\left\{ \left( 1- \Phi\left(\frac{t}{k+1} \right) \right)   \bigwedge \Psi(\frac{x}{k})  \right\}.
\end{equation}
\end{thm}
 \begin{thm}
 Let $R_t$ be a renewal reward process with iid uncertain interarrival times $\xi_1,\xi_2,\ldots  $ and iid uncertain rewards $\eta_1,\eta_2,\ldots $ Assume $(\xi_1,\xi_2,\ldots )$ and $(\eta_1,\eta_2,\ldots )$ are independent uncertain vectors. Then the reward rate
 \begin{equation}
 \frac{R_t}{t}\rightarrow\frac{\eta_1}{\xi_1},
 \end{equation}
 in the sense of convergence in distribution as $t\rightarrow \infty$.
 \end{thm}
 \begin{defn}\cite{liu2013extreme}
 Let $U_t$ be an insurance risk process. Then the ruin index is defined as the uncertain measure that $U_t$ eventually becomes negative, i.e.,
 \begin{equation}
 {\rm Ruin}=\mathscr{M}\left\lbrace \inf_{t\geq 0}U_t < 0 \right\rbrace .
 \end{equation}
 \end{defn}

 \section{The Optimization Model}
Here, we consider the process of variation in the wealth of the primary insurance  as
\begin{equation}\label{Ali1}
U_t^x=u+\left[ x(1+\rho)-(\rho-\theta) \right]  ct-x R_t,
\end{equation}
where, $x \in(0,1]$ is  a  control parameter; and $\rho$ and $\theta$ are reinsurance and insurance safety loads, respectively. It is worthy to remind that term $c$, as the insurance premium, should be chosen in a way to ensure that ruin does not happen almost surely for any amount of the initial capital or reinsurance program.

When a company considers taking a reinsurance contract,  solving the following optimization problem identifies the ruin index.
\begin{equation}
{\rm Ruin}=\max_{k\geq 1}\sup_{z\geq 0}\bigg\{\Phi\Big(\frac{z}{k+1} \Big)\bigwedge  1-\Psi\Big( \frac{u+\big[ x(1+\rho)-(\rho-\theta) \big]c z}{xk}  \Big) \bigg\}.
\end{equation}
Note that solving this problem is very costly; it is almost impossible for most uncertain distributions, and just its approximations would be identified. To overcome this difficulty, we propose a simpler method for the UMR and construct a constraint to prevent ruin.
For the ultimate ruin, the ruin at the infinite time, it is sensible to model an optimization problem dependent on $t$ and taking the objective function as the asymptotic uncertain expected value of the insurer's risk process. Formally, we consider the following constrained optimization problem.
\begin{equation}
\begin{array}{rrclcl}
\displaystyle \max_{0\leq x\leq 1} & \multicolumn{3}{l}{\displaystyle\lim_{t\longrightarrow \infty} \displaystyle\frac{E\left[ U_t^x \right]}{t}  }\\[4mm]
\textrm{s.t.} & {\rm Ruin} \leq \varepsilon,
\end{array}
\end{equation}
where $\varepsilon$ is a small positive number chosen by the administration as maximum tolerable ruin.

\subsection{Calculating the Uncertain Measure of Ruin}\label{Sec1}

First, we need to compute the UMR for primary insurance. In the contract duration, arrival times coincide with some point in the interval $(0,t]$. Since ruin might happen in only  one of these arrivals,  we regard $t=\sum_{i=1}^{N_t}\xi_i$, where $\xi_i$'s are uncertain interarrival times. Therefore, the risk process becomes
\begin{equation}\label{Ali2}
U_t=u+\sum_{i=1}^{N_t}(c\xi_i-\eta_i).
\end{equation}
Let $\xi_i$ and $\eta_i$ be independent uncertain variables. Then by Theorem, the inverse uncertainty distribution of $\eta_1-c\xi_1$ is
\begin{equation}
\Upsilon^{-1}_1(\alpha) = \Psi^{-1}_1(\alpha)-c\Phi^{-1}_1(1-\alpha)=\Psi^{-1}(\alpha)-c\Phi^{-1}(1-\alpha).
\end{equation}
Let $\Theta_i=\eta_i-c\xi_i $ and $M_k = \Theta_1+ \cdots +\Theta_k$.
It can be easily understood that the inverse uncertainty distribution of  $M_k$ for all $k$ is
\begin{eqnarray}
L^{-1}_k(\alpha)=k[\Psi^{-1}(\alpha)-c\Phi^{-1}(1-\alpha)],
\end{eqnarray}
 where $L_k$ is the uncertainty distribution of the uncertain variable $M_k$.

Notice that ruin happens when $ \max_{k\geq 1}\left\lbrace M_k \right\rbrace \geq u.$  Recall that direct calculation of this quantity is a very costly task using its uncertainty distribution, especially when uncertain distributions like lognormal are involved in the model.
Taking $f=\max_k\left\lbrace x_1,\ldots ,x_k\right\rbrace -u$,  $M_k$'s and $f$  satisfy the requirements of Theorem~\ref{Theorem11}. Therefore, the UMR  is
\begin{equation}
\mathscr{M}\left\lbrace f(M_1,M_2,\ldots ,M_k)> 0  \right\rbrace .
\end{equation}
For $k$ big enough, the root $\alpha$ of
\begin{eqnarray}
\max_k \left\lbrace L^{-1}_1(1-\alpha),  L^{-1}_2(1-\alpha),\ldots ,L^{-1}_k(1-\alpha)  \right\rbrace - u=0
\end{eqnarray}
identifies  the UMR. Equivalently, we have to solve
\begin{eqnarray}
\max  \bigg\{ \Psi^{-1}(1-\alpha)-c\Phi^{-1}(\alpha), \ldots,
k(\Psi^{-1}(1-\alpha)-c\Phi^{-1}(\alpha)) \bigg\} =u .
\end{eqnarray}
We need the assumption
\begin{equation}
{\rm Premium} > E[\eta_i]E[\xi_i ],
\end{equation}
in the selection of value  $c$,   because after each interarrival time,   one expects  $E[\xi_i ]$ as the amount of loss, and so the premium must be bigger than this amount.  Observe that  $\eta_i $'s are iid, means that for all $i$ these expectation amounts are identical, as for the $\xi_i $s. This criterion on premiums prevents the model from the obvious event of ruin being determined in a long run.

 To illustrate the practicality of this approach in calculating the ruin, consider the following concrete example.
 \begin{example}\label{example1}
 Consider an insurance company takes some risk in signing a particular class of underwriting. They invite some experts to build two uncertain distributions for the potential number of claims and their severity. Suppose their opinion on the severity of claims and the interarrival times are  $\mathscr{LOGN}(2,1) $  and  $\mathscr{L}(1,3)$, respectively.

 A way to know the measure of possible bankruptcy caused by these contracts is to identify the company's UMR. As it was said, it is the root $\alpha$ of the following equation. To calculate the UMR, we used online Mathematica \cite{Wolfram|One}.
\begin{eqnarray}
\max_{k\geq 1}& &\bigg\{ \exp\Big( 2+\frac{\sqrt{3}}{\pi}\ln\frac{1-\alpha}{\alpha}\Big)-26(2\alpha+1),\ldots, \\ &&\hspace{1cm}k\bigg(\exp\Big( 2+\frac{\sqrt{3}}{\pi}\ln\frac{1-\alpha}{\alpha}\Big)-26 (2\alpha+1)\bigg)\bigg\}=u.
\end{eqnarray}
 Fig.\ref{fig:fig1} represents the UMR concerning the number of claims $k$ with fixed $u=10000$. Notice that this measure gets stable after $k$ gets big enough ( for $k=10000$ and higher), and then the measure does not change considerably. This observation is a completely rational outcome because we expect the ultimate ruin to lift by increasing the number of claims.  After a point on, there is no significant fluctuation in the number UMR, which is mostly because of the premiums' growth. We also illustrated the UMR for different amounts of the initial capital $u$, for fixed $k=100$, Fig.\ref{fig:fig2}. As it is depicting, increasing the initial capital reduces the measure of ruin substantially to a point where it is almost zero.
 \begin{figure}
\begin{subfigure}{.475\textwidth}
  \centering
  \includegraphics[width=0.95\linewidth]{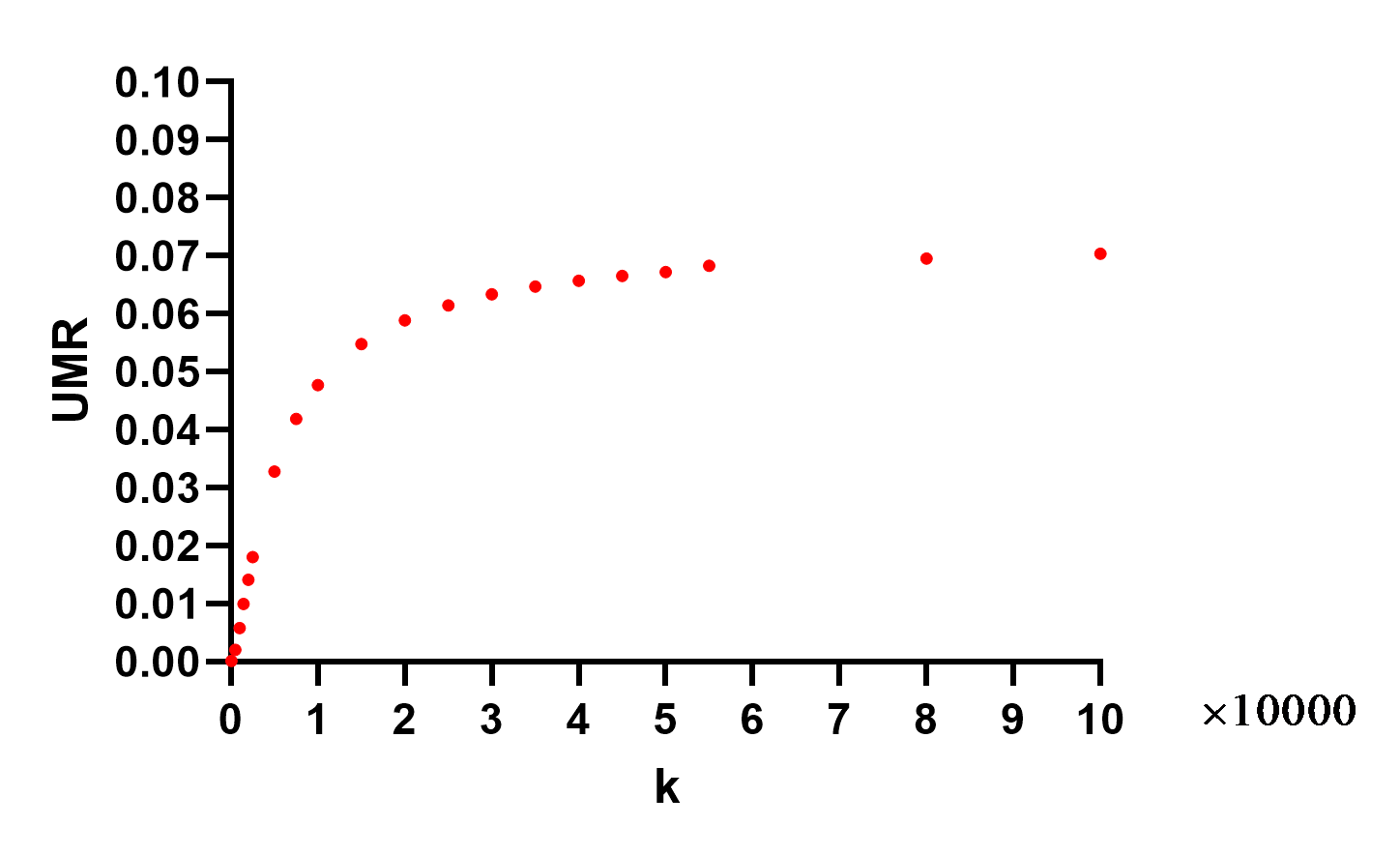}
  \caption{ Number of Claims $k$ for $u=10000$.}
  \label{fig:fig1}
\end{subfigure}%
\hspace{5mm}
\begin{subfigure}{.475\textwidth}
  \centering
  \includegraphics[width=0.95\linewidth]{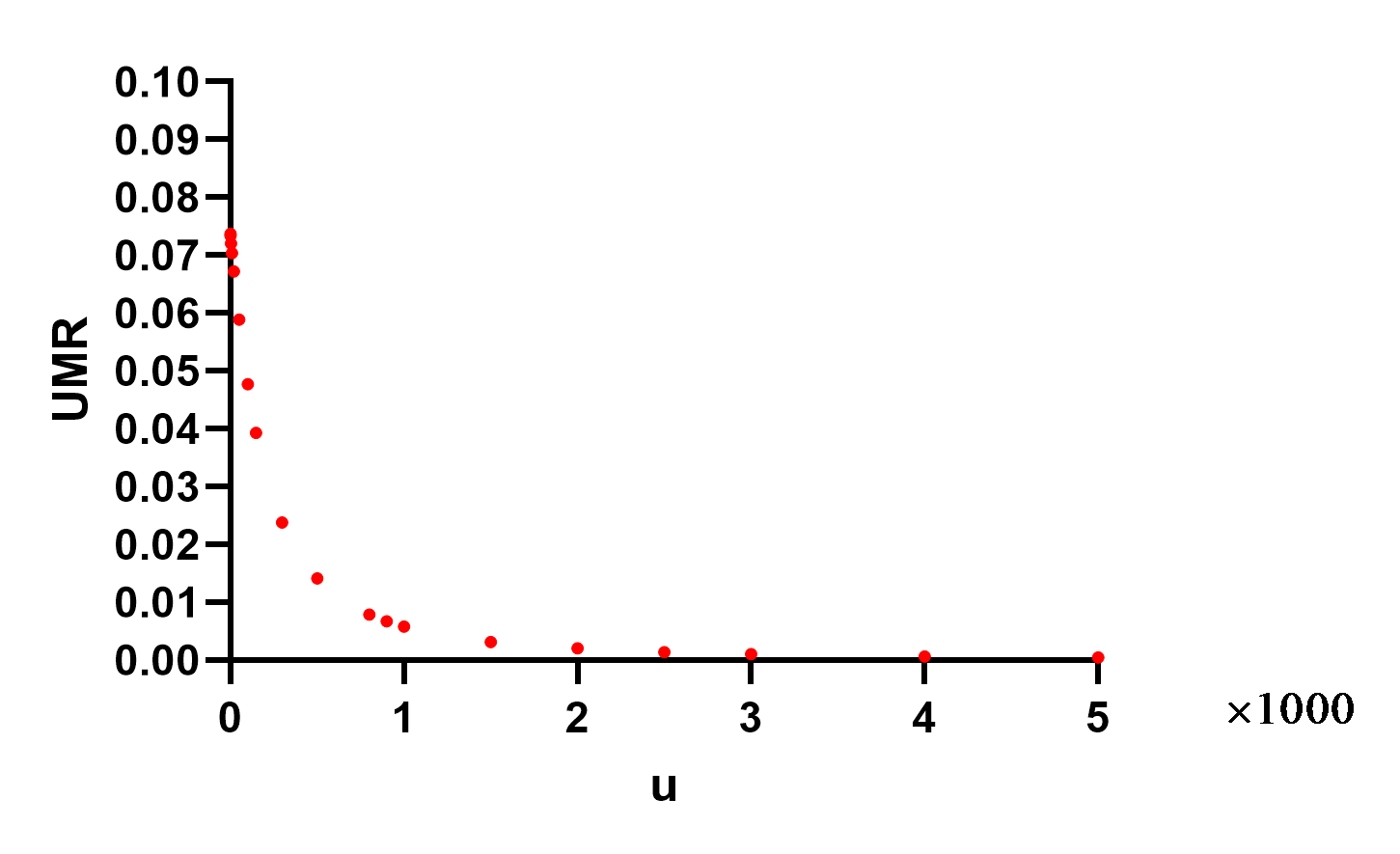}
  \caption{Different amounts of  $u$ for  $k=100$.}
  \label{fig:fig2}
\end{subfigure}
\caption{The sample UMR in Example \ref{example1}.}
\label{fig:fig}
\end{figure}

\end{example}

\subsection{Optimal Retention}
Applying these results for the case of quota share reinsurance contract, we devise a model to explain the primary insurance situation. Solving the optimization problem provides us the optimal retention value $x$. First, we need to calculate the UMR for the case of reinsurance involvement. After sharing its risk with reinsurance, the primary insurance has the following wealth process
\begin{equation}\label{Ali3}
U_t^{x}= u+\sum_{i=1}^{N_t}  \left[ x(1+\rho)-(\rho-\theta) \right] c \xi_i  -x\eta_i,
\end{equation}
with the assumption  $\rho > \theta $, meaning that the risk in reinsurance is more significant than the one in primary insurance. Notice that if we take $x=1$,   the wealth process \eqref{Ali3} describes the primary insurance without reinsurance partnership.

 Let $ \left[ x(1+\rho)-(\rho-\theta) \right]  c$ be denoted by $\beta$. Similar to Sec. \ref{Sec1}, each $x\eta_i-\beta \xi_i$ has the  inverse uncertainty distribution
\begin{equation}
\hat{\Upsilon}^{-1}(\alpha)=x\Psi^{-1}(\alpha)-\beta\Phi^{-1}(1-\alpha).
\end{equation}
Moreover,  for $\hat{\Theta}_i=x\eta_i-\beta\xi_i$, uncertain variables $ \hat{M}_k= \hat{\Theta}_1+ \cdots + \hat{\Theta}_k $
has the following inverse uncertainty distributions
\begin{equation}
\hat{L}_k^{-1}(\alpha)=k\left( x\Psi^{-1}(\alpha)-\beta\Phi^{-1}(1-\alpha)\right),
\end{equation}
where $\hat{L}_k$ is the uncertain distribution of $ \hat{M}_k$.
 To reflect the profitability of the contracts for the primary insurance company, we consider
\begin{equation}
\Bigg[ x(1+\rho)-(\rho-\theta) \Bigg]c-x\int_0^1 \left[ \frac{\Psi^{-1}(\alpha)}{ \Phi^{-1}(1-\alpha)} \right]d\alpha,
\end{equation}
as the objective function,  that is the expected value of the wealth with respect to $t$. It means that the company makes a profit equal to the objective function at each unit of time in the long run. Our purpose is to make this amount positive and as large as possible by an optimal value of $x$. Considering the profitablity of the company in small time intervals influences the idea of optimizing the expected value of the primary insurance company's wealth asymptotically.

The problem identifies the optimal retention rate that maximizes the objective function,  the company's wealth, constrained by keeping its UMR under a predefined value. Thus, we have the following optimization problem.
\begin{equation}
\begin{array}{l}
\displaystyle \max_{0 \leq x \leq 1} \Big[ x(1+\rho)-(\rho-\theta) \Big]c-x\int_0^1  \frac{\Psi^{-1}(\alpha)}{ \Phi^{-1}(1-\alpha)} d\alpha \\[3mm]
\textrm{s.t.}\quad  \textrm{Root}\Big\lbrace \max_{k\geq 1}\big\lbrace \hat{L}^{-1}_1(1-\alpha),\ldots ,\hat{L}^{-1}_k(1-\alpha)  \big\rbrace - u =0  \Big\rbrace\leq \varepsilon.
\end{array}
\end{equation}

\begin{rem}
Observe that in this optimization problem, for $x=0$, the insurance company transfers the whole risk to the reinsurance and leaves with no premium. On the other hand, $x=1$ means that the primary insurance takes the whole risk by itself and does not get into the reinsurance contract. The value of  $x$ is completely dependant on the distributions of severity and number of claims. Recall that an expert provides these uncertainty distributions. While, initial capital, parameters such as $\rho$ and $\theta$ depend on the amount of premium $c$. The value of $c$ is also  provided based on the management policy.
\end{rem}

\begin{example}\label{example2}
We assume an insurance company is deciding on arranging a great amount of non-life insurance contracts. Let it aims to underwrite many houses and businesses in a new area. Since it has no enough experience in this region, the management may ask some experts to contribute their opinions via uncertain distributions on the number and severity of claims based on their experiences.
Consider the case where the interarrival times of the renewal process have a linear distribution, $\mathscr{L}(1,3)$, and claim variables follow the lognormal distribution $\mathscr{LOGN}(2,1)$. This means that experts believe, overall, that for each policyholder there would be at least one claim in the long run, and the number of claims is not bigger than 3. Also, they believe that the amount of loss that each claim produces is in a way that makes a log-normal distribution, meaning that there is always a considerable uncertain measure devoted to rare events.

Now  suppose
$\rho =0.9 , \theta =0.8,  c=26$, and $\varepsilon=0.005$.
With these assumptions,  the optimization problem reads as
\begin{equation}
\begin{array}{l l} 
\displaystyle \max_{0\leq x \leq 1 } & 39.6x-2.6 \\
\textrm{s.t.} & \textrm{Root} \Big\lbrace  \max_{k\geq 1}  k\big[ x \exp(2+\frac{\sqrt{3}}{\pi}\ln\frac{1-\alpha}{\alpha })\\
& \hspace{3cm}-(49.4x-2.6)(2\alpha+1) \big] -u \Big\rbrace \leq 0.005.
\end{array}
\end{equation}
The optimal retentions of the company for $k=100$, a fixed number of possible claims, and for different initial capital amounts are represented in Fig.\ref{figg:fig3}, which show the behavior of retention rate with varying this value. To analysis the effect of claims' number on the retention rate, we fix the initial capital at  $u=100000$.

When $u$ is small the UMR is significantly high, consequently, the retention needs to be small to keep the UMR under the assumed limit. Raising the initial capital certainly gives a more chance to the insurance company to take more risks and at the same time keeping the UMR under the limit, Graph (\ref{figg:fig3}). Investigating Graph (\ref{figg:fig4}) demonstrates that, growing the number $k$, increases the UMR, which makes the retention rate to decrease in a way that the problem's constraint stays satisfied. Almost by $k=7000$ (the fixed initial capital), and the numbers higher than that, the amount of retention gets stabilized.

\begin{figure}
\begin{subfigure}{.475\textwidth}
  \centering
  \includegraphics[width=0.95\linewidth]{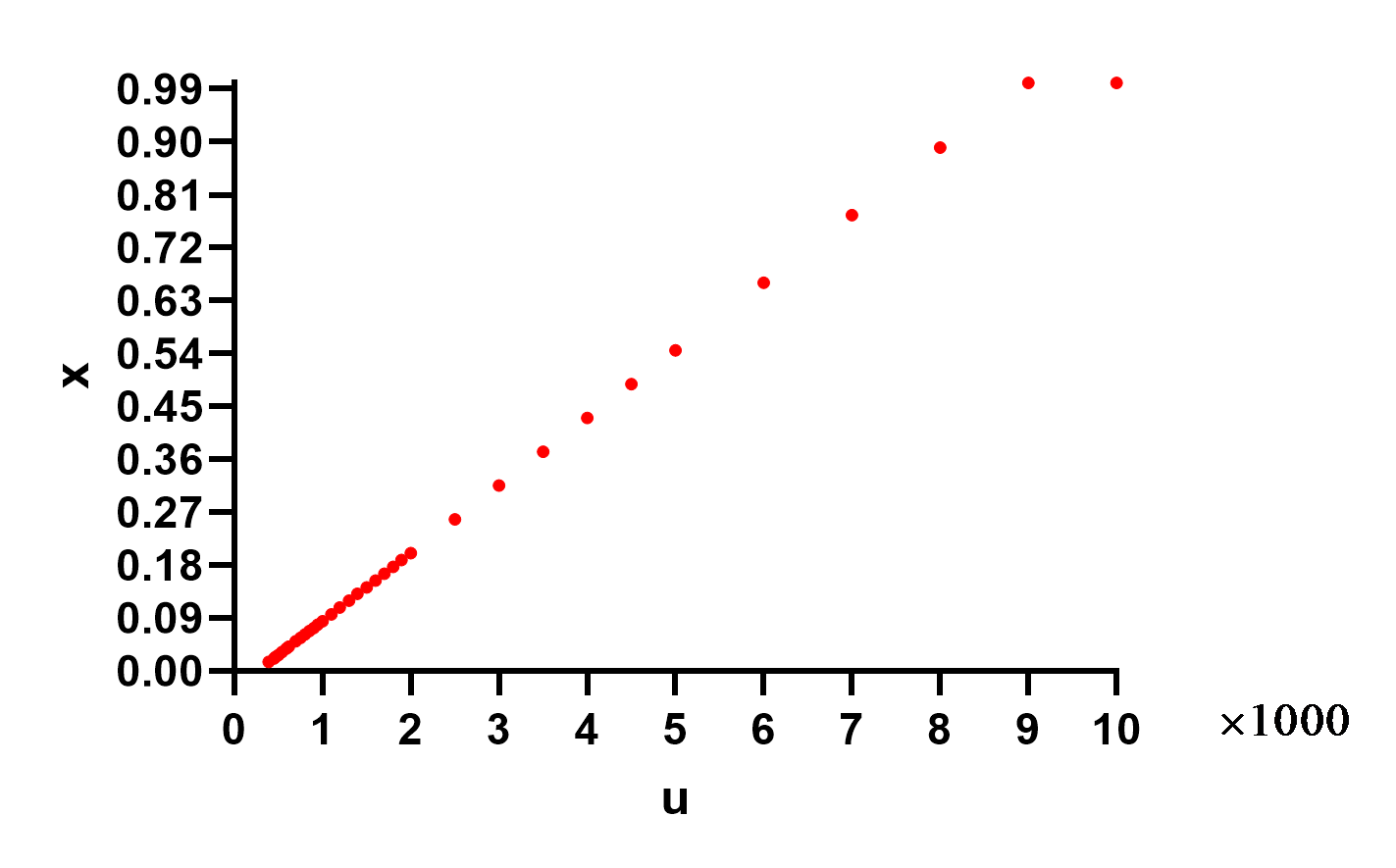}
  \caption{Different Amount of the Initial Capital.}
  \label{figg:fig3}
\end{subfigure}%
\hspace{5mm}
\begin{subfigure}{.475\textwidth}
  \centering
  \includegraphics[width=0.95\linewidth]{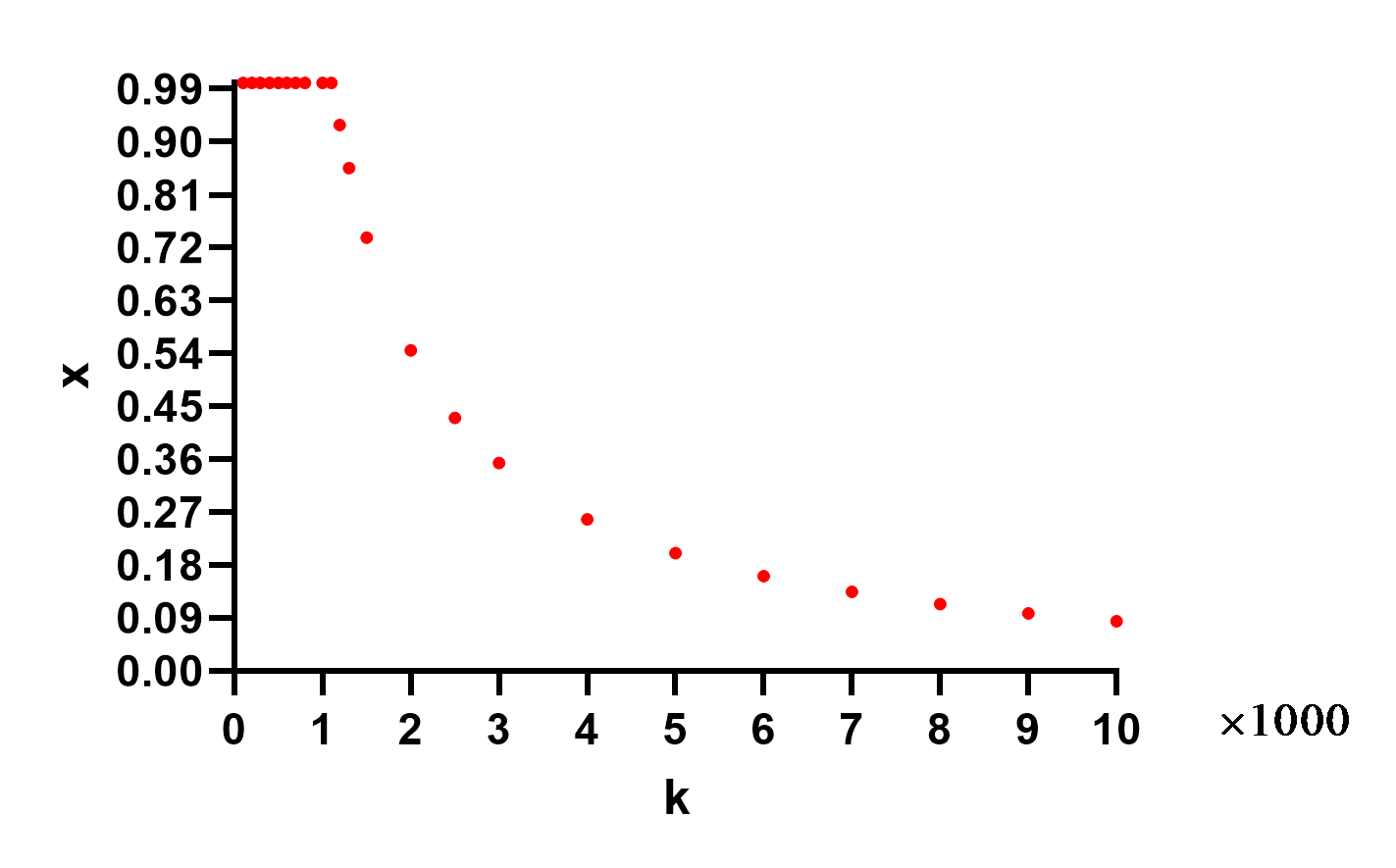}
  \caption{ Different Possible  Amount of the Claims. }
  \label{figg:fig4}
\end{subfigure}
\caption{Optimal Retention Rate for the Insurance Company Constrained to Keep UMR Under Preassumed Number in Example \ref{example2}.}
\label{figg:figg}
\end{figure}

\end{example}

 Notice that the number of sold contracts by the primary insurance company can be considered as the potential claims' number. With this assumption, one can assume that a catastrophe situation would occur when all these contracts have a loss and confirm a claim.

This optimization problem can also be more practically used by incorporating other optimal control parameters to design contracts and reduce the risk while raising their potential gain. It should also be noticed that increasing the initial capital $u$  lessens the risk of ruin, too. However, we don't use it as a control parameter as it is absent in the objective function.

\subsection{A Perturbed Extension}
Following the idea of perturbing the risk process by a Brownian motion~\cite{dufresne1991risk}, one can perturb the UMR by the  Liu process. By perturbation, we mean adding an uncertain process that can take positive and negative numbers. This generalized model is used to include more practical considerations in the model. We can interpret it as the imperfect human reckoning in the expert opinion's contribution when only a few data is available. In this way, adding a potential fluctuation to the model allows us to mimic this incompetence. It can also be recognized as the variation caused by dependency between uncertain variables used in the independent model, recalling that they are mostly dependent in real life.

The perturbed problem would be modeled  as
\begin{equation}
U_t=u+ct -R_t+C_t,
\end{equation}
where $C_t$ is a Liu process having the following inverse uncertainty distribution.
\begin{equation}
\tilde{\Phi}_t^{-1}(\alpha)=t\frac{\sqrt{3}}{\pi}\ln(\frac{\alpha}{1-\alpha}).
\end{equation}
We assume that the Liu process, at each point of time, is independent of other uncertain variables in the model. Following the inference for the non-perturbed case, we consider the uncertain variables
\begin{equation}
\tilde{\Theta}_i=\eta_i-c\xi_i - C_i ,
\end{equation}
as a  perturbed version of $\Theta$. The inverse uncertainty distribution of this uncertain variable  is
\begin{equation}
\tilde{\Upsilon}^{-1}_i(\alpha) = \Psi^{-1}_i(\alpha)-c\Phi^{-1}_i(1-\alpha)-\tilde{\Phi}_i^{-1}(1-\alpha).
\end{equation}
Defining
$ \tilde{M}_k=\tilde{\Theta}_1 +\tilde{\Theta}_2+\cdots+\tilde{\Theta}_k, $
the UMR becomes the uncertain measure of the following event
\begin{equation}
\max_{k\geq 1}\left\lbrace \tilde{M}_k \right\rbrace \geq u,
\end{equation}
with the inverse uncertainty distribution
\begin{eqnarray}
\tilde{L}^{-1}_k(\alpha)&=&k(\Psi^{-1}(\alpha)-c\Phi^{-1}(1-\alpha)-\tilde{\Phi}_k^{-1}(1-\alpha)), \qquad k\geq 1.
\end{eqnarray}
 Using Theorem \ref{Theorem11}, for $k$ big enough, the UMR is the root $\alpha$ of
 \begin{eqnarray}
\max_k \left\lbrace \tilde{L}^{-1}_1(1-\alpha),  \tilde{L}^{-1}_2(1-\alpha),\ldots ,\tilde{L}^{-1}_k(1-\alpha)  \right\rbrace - u =0,
\end{eqnarray}
or equivalently the root of
\begin{equation}
\begin{array}{l}
\max_{k\geq 1}  \Big \{ \Psi^{-1}(1-\alpha)-c\Phi^{-1}(\alpha)-\tilde{\Phi}_k^{-1}(\alpha),
\ldots, \\ \hspace{1.5cm} k(\Psi^{-1}(1-\alpha)-c\Phi^{-1}(\alpha)-\tilde{\Phi}_k^{-1}(\alpha)) \Big\} =u .
\end{array}
\end{equation}
Similarly, when  reinsurance is involved in the model, we have
\begin{equation}
U_t^x=u+\left[ x(1+\rho)-(\rho-\theta) \right]  ct-x R_t+C_t .
\end{equation}
 Defining $\tilde{\Theta}_i=x\eta_i-\beta\xi_i-C_i$, and
$ \tilde{M}_k = \tilde{\Theta}_1+\tilde{\Theta}_2 + \cdots + \tilde{\Theta}_k, $
with the inverse uncertainty distribution
\begin{equation}
\tilde{L}_k^{-1}(\alpha)=k\left( x\Psi^{-1}(\alpha)-\beta\Phi^{-1}(1-\alpha)-\tilde{\Phi}_k^{-1}(1-\alpha))\right),
\end{equation}
the UMR  is the root $\alpha$ of
\begin{equation}
\max_{k\geq 1} \left\lbrace \tilde{L}^{-1}_1(1-\alpha),  \tilde{L}^{-1}_2(1-\alpha),\ldots ,\tilde{L}^{-1}_k(1-\alpha)  \right\rbrace - u =0 .
\end{equation}
\section{Conclusions}\label{conclusions}
In this paper, we investigated a model designed for an insurance company that wants to share a part of its risk with a reinsurer. Uncertainty theory was applied, for the first time, to model the involved uncertainty of the problem. Some comparisons between initial capital, the maximum number of possible claims, and the optimal retention were achieved. Further research can model other kinds of reinsurance contracts as an uncertain optimization problem with similar or different sorts of constraints.


%
 \section*{Conflict of interest}
 The authors declare that they have no conflict of interest.



\end{document}